\documentclass[ twoside, 11pt]{article}
\usepackage[margin=1in]{geometry}

\usepackage{setspace}
\usepackage{amsmath}
\usepackage{amssymb}
\usepackage[square]{natbib}
\setcitestyle{numbers}

\newtheorem{theorem}{Theorem}[section]
\newtheorem{corollary}{Corollary}[section]
\newtheorem{lemma}{Lemma}[section]

\newtheorem{example}{Example}[section]

\newcommand{\Proof}{\textbf{Proof. }}            
\newcommand{\Proofof}[1]{\textbf{Proof of #1. }} 

\newcommand{\RR}{\mathbb{R}}

\newcommand{\ZZ}{\mathbb{Z}}

\newcommand{\ii}{{\mathrm i}}
\newcommand{\ee}{{\mathrm e}}

\newcommand{\dd}{{\mathrm{d}}}

\newcommand{\BI}{\mathrm{Bi}}    
\newcommand{\NB}{\mathrm{NB}}   

\newcommand{\norm}[1]{\|#1\|}                    
\newcommand{\ab}[1]{\vert#1\vert}                
\newcommand{\Ab}[1]{\Big\vert#1\Big\vert}        

\newcommand{\exponent}[1]{\exp\{#1\}}            
\newcommand{\Exponent}[1]{\exp\Bigl\{#1\Bigr\}}  

\newcommand{\Expect}{\mathrm{E}}                 
\newcommand{\Prob}{\mathrm{P}}                   

\newcommand{\qubar}{\overline{q}}
\newcommand{\pbar}{\overline{p}}

\newcommand{\vfi}{\varphi}
\newcommand{\wE}{\w\Expect}

\newcommand{\w}{\widehat}

\newcommand{\floor}[1]{{\lfloor #1\rfloor}}

\begin{document}
\title{On large deviations for sums of discrete m-dependent
random variables}

\author{V. \v Cekanavi\v cius  and P. Vellaisamy      \\
{\small
Department of Mathematics and Informatics, Vilnius University,}\\
{\small Naugarduko 24, Vilnius 03225, Lithuania.}\\{\small E-mail:
vydas.cekanavicius@mif.vu.lt } \\{\small and}
\\{\small
 Department of Mathematics, Indian Institute of Technology Bombay,} \\
 {\small Powai, Mumbai-
400076, India.}\\{\small  E-mail: pv@math.iitb.ac.in} }
\date{}

\maketitle

\begin{abstract}

  The ratio  $P(S_n=x)/P(Z_n=x)$ is investigated for three cases: (a) when $S_n$ is a  sum of 1-dependent non-negative integer-valued random variables (rvs), satisfying some moment conditions, and $Z_n$ is Poisson rv; (b) when $S_n$ is a statistic of 2-runs and $Z_n$ is  negative binomial rv; and (c) when $S_n$ is statistic of $N(1,1)$-events and $Z_n$ is a binomial r.v. We also consider the  approximation of  $P(S_n\geqslant x)$ by Poisson distribution with parameter depending  on $x$.

\vspace*{.5cm} \noindent {\emph{Key words:} \small
 binomial distribution, large deviations,  m-dependent random variables, negative binomial distribution, Poisson distribution.}

\vspace*{.5cm} \noindent {\small {\it MSC 2000 Subject
Classification}:
Primary 60F10.   
Secondary 60G50;     
}
\end{abstract}

\newpage

\section{Introduction}

In this paper, we prove some large deviation results for the sum $S_n= X_1+X_2+\cdots+X_n$ of
identically distributed 1-dependent random variables (rvs)
concentrated on non-negative integers.
Recall that a sequence of random variables
$\{X_k \}_{k \geq 1}$
  is called $m$-dependent if, for $1 < s < t < \infty$, $t- s > m$, the sigma-algebras
generated by $X_1,\dots,X_s$ and $X_t,X_{t+1}\dots$ are
independent. It is clear that,  by grouping consecutive summands, we
can reduce the sum of $m$-dependent variables to the sum of
1-dependent ones.
There are numerous large deviation (LD) results for normal approximation, see, for example,  an overview of the
classical results in \cite{Na79} or Chapter VIII in \cite{Pet75}. The
LD results for sums of $m$-dependent rvs, when approximating rv is normal,  have been considered in \cite{Gh74, H82} and Theorems 4.29--4.31 in \cite{SAST91}. The LD for Poisson distribution are  not so comprehensively studied. For the case of independent integer-valued variables, see \cite{ALST03, CeVa99, ChCh92, De92, Iv75,  WoMi83} and the references therein.   The LD problems for Poisson approximation to the tail of sums of dependent Bernoulli variables via Stein's method is investigated in  \cite{CFS13}.

 Apart from the LD results for tails of distributions,  local versions of LD (the so-called Richter's type LD) exist for the
 ratio of two densities, or the ratio of probability and density, see, for example, Chapter 7 in \cite{IL71}.
 For two lattice distributions, it is more natural  to compare their probabilities.  In this paper, we compare $P(S_n=x)$ to the probability of approximating a discrete distribution, assuming $\Expect S_n$ to be  large and $\Expect X_i$ to be small. Unlike in \cite{CFS13}, we consider rvs which are not necessarily Bernoulli rvs.

\section{Results}

We use similar assumptions to the ones from \cite{CeVe15a}, where the closeness of $S_n$ to various discrete distributions has been estimated. Further on, all absolute positive
constants are denoted by the same symbol $C$. Sometimes to avoid possible ambiguities, the constants $C$ are supplied
with indices. Notation $\theta$ is used for any real or complex number satisfying $\ab{\theta}\leqslant 1$.

 Let $X_1,X_2,\dots,X_n$ be  identically distributed   non-negative  integer-valued 1-dependent rvs. The first tow factorial moments of $X_1$ are respectively denoted by
$\nu_1:=\Expect X_1$ and $\nu_2:=\Expect X_1(X_1-1)$.
 Let $y:=(x-n\nu_1)/(n\nu_1)$,
 \begin{equation}\label{neqn1}
 \Lambda(y):=-n(1-\nu_1)\sum_{j=2}^\infty\bigg(\frac{\nu_1}{1-\nu_1}\bigg)^j\frac{y^j}{j(j-1)},\; \;
 \gamma:=\ee^{(1.5)C_0}\max(\nu_1^2,\nu_2,\Expect X_1X_2)
 \end{equation}
 and let $Z_n$ denote the  Poisson rv with parameter $n\nu_1$.
 We first make the following assumptions.  For some absolute constant $C_0\geqslant 1$ and positive integer $x$,
 \begin{eqnarray}
 &&X_1\leqslant C_0,\quad \ee^{5C_0}\nu_1\leqslant 0.002,
 \quad \ab{y}\leqslant \frac{1}{10}, \label{sal1}\\
 && \nu_2\leqslant\frac{\nu_1\ee^{-(1.5)C_0}}{20},\quad \Expect X_1X_2\leqslant\frac{\nu_1\ee^{-(1.5)C_0}}{20}\label{sal2}.
 \end{eqnarray}
 The smallness of $\nu_1$ is determined by the method of proof. Note that $\nu_1$, though very small,  can nevertheless be of constant order. 


 \begin{theorem}\label{THM1} Let the assumptions in (\ref{sal1})--(\ref{sal2}) hold. Then
 \begin{equation}
 \frac{\Prob(S_n=x)}{\Prob(Z_n=x)}=\ee^{\Lambda(y)}(1+C_1\gamma y^2)^n\Big(1+\theta C_2\gamma(\nu_1^{-1}+ny^2)\Big).
  \label{t1}\end{equation}
 \end{theorem}

  Note the estimate in (\ref{t1})  is similar in form to the estimates of \cite{ALST03} and \cite{Iv75}, since  it contains some equivalent of Cramer's series in the exponent.


  \begin{corollary} \label{cor1} Let, in Theorem \ref{THM1},  $x-n\nu_1=o(\sqrt{n\nu_1})$, $\gamma=o(\nu_1)$, $\nu_1=o(1)$ and $n\nu_1\to\infty$, as $n\to\infty$. Then
    \begin{equation}
 \frac{\Prob(S_n=x)}{\Prob(Z_n=x)}=1+o(1).
   \label{t2}\end{equation}
    \end{corollary}


 \begin{example} \label{exmpl}
Consider a sequence of independent Bernoulli trials. The event $(k_1,k_2)$ occurs if $k_1$ consecutive failures are
followed by $k_2$ consecutive successes. More formally, let $\eta_i$
be independent Bernoulli
 $Be(p)$ ($0<p<1$)  variables and
$U_j=(1-\eta_{j-m+1})\cdots
(1-\eta_{j-k_2})\eta_{j-k_2+1}\cdots\eta_{j-1}\eta_j$,
$j=m,m+1,\dots,n$, where $m=k_1+k_2$ and $k_1>0$
and $k_2>0$ are fixed integers. Then
$N(n;k_1,k_2)=U_m+U_{m+1}+\cdots+U_n$  denotes the number of $(k_1, k_2)$ events in
 $n$ Bernoulli trials.
 %
 %
%
%

It is well known that $N(n;k_1,k_2)$
has limiting Poisson distribution, see \cite{HuTs91, Vel04}. Note also that $U_1,U_2,\dots$ are
$m$-dependent. Consequently, if $m>1$, then Theorem \ref{THM1} can
not be applied directly. However, if we consider a subsequence  $n=(K+1)m-1$, then the summands can be grouped in
the following natural way:
\begin{eqnarray*}
N(n;k_1,k_2)&=&(U_m+U_{m+1}+\cdots+U_{2m-1})+(U_{2m}+U_{2m+1}+\cdots+U_{3m-1})+\dots \\
             &=&X_1+X_2+\dots+X_K.
\end{eqnarray*}
\noindent Here, each $X_j$ contains $m$
 summands and $X_1,X_2,\dots,X_K$ are 1-dependent.
  Let  $\alpha(p)=(1-p)^{k_1}p^{k_2}$. Then, for $j=1,\dots, K$
\begin{equation*}
X_j= \begin{cases} 1, & \mbox{with probability } m\alpha(p), \\
0, & \mbox{with probability } 1-m\alpha(p), \end{cases}
\end{equation*}
 $\nu_2(j)=\nu_2(K+1)=0$,
$\nu_1(j)=m\alpha(p)$
and  $\Expect X_1X_2=\alpha(p)^2m(m+1)/2$, see \cite{CeVe15b}.

Let $\alpha(p)\to 0$,  $n\to\infty$. Then, for large $K$, conditions (\ref{sal1})--(\ref{sal2}) are satisfied  and, consequently, Theorem \ref{THM1} can be applied. In particular, if  $p\to 0$, $np\to\infty$, then the equivalence relation (\ref{t2}) holds for $N(1,1)$ in the interval $\ab{x-np(1-p)}<C\sqrt{np}$.

\end{example}


An extension of Theorem \ref{THM1} for the tails of distributions by our method of proof is somewhat complicated, as  we are unaware about a satisfactory analogue of Mill's ratio for Poisson distribution. Therefore, we restrict ourselves just to one result, which is of interest, mainly due to an  atypical choice of parameters of 
 approximation. In all standard LD results, for all values of argument $x$, the same approximating distribution is used.
Such an approach is natural for  Kolmogorov, total variation or similar metrics.
 For the LD,   a choice of approximation
 is not so obvious.  From a practical point of view, it make sense to use $\Prob(Z^*\geqslant x)$, with $Z^*=Z^*(x)$.
 For independent rvs, the LD results  with Poisson approximations depending on $x$ are obtained in \cite{BaJe89} and \cite{De92}.  We use similar approach for $m$-dependent rvs. Let $Z_n^{*}$ be the Poisson variable with parameter $n\lambda^{*}$ and let
\[\lambda^{*}:=\nu_1\bigg(1-\frac{\nu_1}{1-\nu_1}y\bigg),\qquad
\Lambda^{*}(y):=n\sum_{j=2}^\infty\frac{1}{j}\bigg(\frac{\nu_1}{1-\nu_1}\bigg)^jy^j.\]
We first consider the right-tail probability of $S_n$.


\begin{theorem}\label{THM2} Let $x>n\nu_1$ and let the assumptions in (\ref{sal1})--(\ref{sal2}) be satisfied. Then
\begin{equation*}
\frac{\Prob(S_n\geqslant x)}{\Prob(Z^{*}_n\geqslant x)}=\ee^{\Lambda^{*}(y)}(1+C_3\gamma y^2)^n\Big(1+C_4\theta\gamma\sqrt{n\nu_1}(\nu_1^{-1}+ny^2)\Big).
\end{equation*}
\end{theorem}

\begin{corollary} \label{cor2} Let, in Theorem \ref{THM2}, $x-n\nu_1=O(\sqrt{n\nu_1})$, $\gamma=O(\nu_1^2)$ and  $n\nu_1^3=o(1)$, as $n\to\infty$. Then
    \begin{equation}
 \frac{\Prob(S_n\geqslant x)}{\Prob(Z^{*}_n\geqslant x)}=1+o(1).
   \label{t3}\end{equation}
    \end{corollary}

When we deal with  closeness of distributions in total variation, the Poisson approximation is not always the best choice, see, for example, \cite{CeVe15a}. Considering two statistics with an explicitly defined dependency of variables, we will demonstrate that the same is true for LD. First we investigate 2-runs statistic, which is one of the most popular examples of sums of 1-dependent variables, see, for example, \cite{CFS13, PeCe10, WXia08} and the references therein.

 Let
$S_\xi=\xi_1+\xi_2+\dots+\xi_n$, where $\xi_i=\eta_i\eta_{i+1}$
and $\eta_i\sim Be(p), 1 \leq i \leq n+1,$ are iid
Bernoulli variables with parameter $p=\Prob(\eta_1=1)=1-\Prob(\eta_1=0)$.
In \cite{CFS13}, the tail of $S_\xi$ was compared to the tail of  Poisson distribution.
 Meanwhile, it is known that, for total variation, $S_\xi$ is much closer to the negative binomial rv, see  \cite{WXia08}. We investigate the benefits that can be gained, if Poisson approximation in LD is replaced by the negative binomial one. Let the negative binomial distribution and its characteristic function  be defined by
\begin{equation}\label{nbinom}
\NB(r,\bar{q})\{j\}=\frac{\Gamma(r+j)}{j!\Gamma(r)}\,\bar{q}^r\pbar^j,
\qquad
\w{\NB}(r,\qubar)(\ii t)=\bigg(\frac{\qubar}{1-\pbar\ee^{\ii t}}\bigg)^r,
\end{equation}
respectively. Here $\qubar+\pbar=1$, $j=1,2,\dots$.
The choice of parameters  that ensures matching of the means and variances of $S_\xi$ and its approximation is:
\begin{equation}\label{NBpar}
\pbar=\frac{np^3(2-3p)-2p^3(1-p)}{np^2+np^3(2-3p)-2p^3(1-p)},\quad r=\frac{n^2p^4}{np^3(2-3p)-2p^3(1-p)}.
\end{equation}

\begin{theorem}\label{THM3} Let $p\to 0$, $np\to\infty$, $\ab{x-np^2}=o(\min(np^2,n^{2/3}p^{2/3}))$,
as $n\to\infty$. Then
\[
\frac{\Prob(S_\xi=x)}{\NB(r,\bar{q})\{x\}}=1+o(1). \]
\end{theorem}

 Observe that the zone of equivalence for probabilities in Theorem \ref{THM3} is larger than the one in Theorem \ref{THM1}. For example, if $p=n^{-1/4}$ then we have $\ab{x-np^2}=o(n^{1/2})$ in Theorem \ref{THM3}, whereas $\ab{x-np^2}=o(n^{1/4})$ in Theorem \ref{THM1}.

 To the best of our knowledge, Theorem \ref{THM3} is the first LD result for the negative binomial distribution, at least for 1-dependent variables.
The larger equivalence zone in Theorem \ref{THM3} is the result of a better matching of corresponding moments.
We have used the negative binomial approximation, which takes care of two matching moments, while the Poisson approximation matches only the mean.
For independent summands, the same effect in LD has been achieved, when the Poisson approximation was replaced by a  shifted Poisson law, see \cite{CeVa99}.

In our last example, we have demonstrated that binomial distribution 
also can outperform Poisson approximation, if the approximate matching of two moments is ensured. We consider statistic of $N(1,1)$-events as defined in Example
\ref{exmpl}. For the sake of simplicity we take $n+1$ iid  Bernoulli variable $\eta_1,\dots,\eta_{n+1}$, $\Prob(\eta_1=1)=p=1-\Prob(\eta_1=0)$
and set $\tilde S=\sum_{j=1}^{n}\eta_j(1-\eta_{j+1})$. In our case, $\alpha=\alpha(p)=p(1-p)$.

We define
the binomial distribution of this paper as
\begin{equation*}\label{binom}
\BI(N,\tilde{p})\{k\}={{N}\choose{ k}}\tilde p^{k}(1-\tilde p)^{N-k},
 \quad
 N=\floor{\tilde N}, \quad \tilde N
=\frac{n^2}{3n-2},
\quad\tilde{p}=\frac{n\alpha}{N}.
\end{equation*}
 Here, we use
$\floor{\tilde N}$ to denote the largest integer less than or equal to  $\tilde N$. The characteristic function of $\BI(N,\tilde{p})$ is equal to
\begin{equation}
\w{\BI}(N,\tilde p)(\ii t)=\big(1-\tilde{p}+\tilde{p}\ee^{\ii t}\big)^N.
\label{chf}
\end{equation}
Note that $N$ must be an integer and therefore the choice of parameters ensures the exact matching of the means and only approximate matching of the variances.

\begin{theorem}\label{THM4} Let $p\to 0$, $np\to\infty$ and $\ab{x-n\alpha}=o(\min(np,n^{2/3}))$, as $n\to\infty$. Then
\[\frac{\Prob(\tilde S=x)}{\BI(N,\tilde p)\{x\}}=1+o(1).\]
\end{theorem}

If $p=n^{-1/3}$, the equivalence zone is $\ab{x-n\alpha}=o(n^{2/3})$ which is much larger than the zone $\ab{x-n\alpha}=o(n^{1/3})$ for Poisson approximation given in $(\ref{t1})$.

\section{Auxiliary results}

Henceforth, let $\w F_n(u):=\Expect \ee^{uS_n}$. Observe that $\w F_n(\ii t)$ is the characteristic function of $S_n$.
In all the formulas, $h$ denotes a solution to the saddle point equation for the binomial distribution:
$(\ln(1-\nu_1+\nu_1\ee^h)^n)'=x$. It is easy to check that
\begin{equation}\label{au1}
x=\frac{n\nu_1\ee^h}{1+\nu_1(\ee^h-1)},\quad \ee^h-1=\frac{y}{1-\nu_1-\nu_1y},\quad\frac{1+\nu_1(\ee^{\ii t+h}-1)}{1+\nu_1(\ee^{h}-1)}=1+\frac{x}{n}(\ee^{\ii t}-1),
\end{equation}
where, as before, $y=(x-n\nu_1)/(n\nu_1)$.
If the assumptions in (\ref{sal1})-(\ref{sal2}) hold, then
$h\leqslant\ee^h-1\leqslant 0.1003$.  Let $u:=h+\ii t$. Then
\begin{eqnarray}
\ab{\ee^u-1}&\leqslant&\ab{\ee^{\ii t}(\ee^h-1)+(\ee^{\ii t}-1)}\leqslant \ee^h-1+\ab{\ee^{\ii t}-1}\leqslant
1.003y+2\ab{\sin\frac{t}{2}},\label{au3}\\
\ab{\ee^u-1}^2&\leqslant&2(\ee^h-1)^2+2\ab{\ee^{\ii t}-1}^2\leqslant 2.02y^2+8\sin^2(t/2). \label{au4}
\end{eqnarray}

We  apply Heinrich's version of the characteristic function method and conjugate distribution, see \cite{H82, H87}. First, we introduce the necessary notations. Let $\{Y_{k}\}_{k \geq 1}$
be a sequence of arbitrary real or complex-valued random
variables. We assume that $\w\Expect (Y_1)=\Expect Y_1$ and, for
$k\geqslant 2$, define $\w\Expect (Y_1,Y_2,\cdots Y_k)$  by
\begin{equation}
 \w\Expect (Y_1,Y_2,\cdots, Y_k)=\Expect Y_1Y_2\cdots
Y_k-\sum_{j=1}^{k-1}\w\Expect (Y_1,\cdots ,Y_j)\Expect
Y_{j+1}\cdots Y_{k}. \label{capY}
\end{equation}

The next lemma is proved in \cite{H82}.
\begin{lemma} \label{Hei3aa} Let $Z_1,Z_2,\dots,Z_k$ be 1-dependent complex-valued random variables with
$\Expect\ab{Y_m}^2<\infty$,  $1 \leqslant m \leqslant k. $ Then
\[
\ab{\wE (Y_1, Y_2, \cdots, Y_k)}\leqslant
2^{k-1}\prod_{m=1}^k(\Expect\ab{Y_m}^2)^{1/2}.
\]
\end{lemma}

In the next lemma, we have collected some facts about  $Y_j:=\ee^{uX_j}-1$. All the derivatives in Lemma \ref{lem1} are with respect to $t$.
\begin{lemma}\label{lem1} Let the assumptions in (\ref{sal1})-(\ref{sal2}) hold. Then, for any $k=1,2,\dots,n$,
\begin{eqnarray}
\ab{Y_k}&\leqslant&\ee^{0.1003C_0}X_k(2\ab{\sin(t/2)}+1.003y),\label{lau1}\\
\quad \ab{Y_k}&\leqslant& \ee^{0.1003C_0}2.1003X_k,\label{lau1aa}\\
\ab{Y_k'}&\leqslant& CX_k,\quad\Expect\ab{Y'_k}\leqslant C\nu_1,\qquad \Expect\ab{Y'_k}^2\leqslant C\nu_1,\label{lau1a}
\\
\Expect\ab{Y_k}^2&\leqslant&\ee^{1.21C_0}\nu_1(8\sin^2(t/2)+2.02y^2), \ \Expect\ab{Y_k}^2\leqslant\ee^{-4.7C_0}0.009,\label{lau2}
\\
\Expect Y_k&=&\nu_1(\ee^{u}-1)+\theta 0.1\nu_1\sin^2(t/2)+\theta 1.01\nu_2\ee^{0.2C_0}y^2,\label{lau3}
\\
\Expect Y'_k&=&\ii \ee^u\nu_1+\theta C\nu_2(\ab{\sin(t/2)}+y),\label{lau3a}
\\
\ab{\wE (Y_1,Y_2)}&\leqslant& 0.12\nu_1\sin^2(t/2)+2.5\ee^{0.3C_0}\Expect X_1X_2y^2,\label{lau4}\\
\ab{(\wE (Y_1,Y_2))'}&\leqslant&C\gamma(\ab{\sin(t/2)}+y),\label{lau4a}\\
\ab{\wE (Y_1,Y_2,Y_3)}&\leqslant&0.7\nu_1\sin^2(t/2)+6\ee^{1.4C_0}\Expect X_1X_2y^2+0.0001\nu_1^2y^2,\label{lau5}\\
\ab{(\wE (Y_1,Y_2,Y_3))'}&\leqslant&C\gamma(\sin^2(t/2)+y^2),\label{lau5a}\\
\ab{\wE (Y_j,  \cdots, Y_k)}&\leqslant& C\gamma(\sin^2(t/2)+y^2)(0.02)^{k-j},\quad j-k\geqslant 1,\label{lau6a}\\
\ab{\wE (Y_j,  \cdots, Y_k)}&\leqslant& 0.16\nu_1\sin^2(t/2)(0.02)^{k-j-3}+65.3\nu_1^2y^2(0.061)^{k-j-3},
\label{lau6}\\
\ab{(\wE (Y_j, \cdots, Y_k))'}&\leqslant&C\nu_1^2(\ab{\sin(t/2)}^3+y^3)(k-j+1)(0.03)^{k-j},\ j-k\geqslant 3.\label{lau7}
\end{eqnarray}
\end{lemma}
\emph{Proof.} The estimate in (\ref{lau2}) follows from (\ref{lau1}). To prove the estimates in (\ref{lau1}), observe that
\[
\ab{Y_k}\leqslant\ab{\ee^{hX_k}(\ee^{\ii t X_k}-1)}+(\ee^{hX_k}-1)\leqslant \ee^{hX_k}X_k(\ab{\ee^{\ii t}-1}+(\ee^h-1))\]
and
\[\ab{Y_k}\leqslant \ee^{hX_k}\ab{\ee^{\ii t X_k}-1}+\ab{\ee^{hX_k}-1}\leqslant \ee^{hC_0}X_k\ab{\ee^{\ii t}-1}+hX_k\ee^{hX_k},\]
and it remains to apply (\ref{au3}). The estimates (\ref{lau1a}) are proved in \cite{CeVe15b}, Lemma 4.3.
The  estimates (\ref{lau3}) and (\ref{lau3a}) follow from the expansion in factorial moments, namely,
\begin{equation}\label{ho}\Expect Y_i=\nu_1(\ee^u-1)+\theta\ee^{hC_0}\nu_2\frac{\ab{\ee^u-1}^2}{2},\end{equation}
which is proved in \cite{CeVe15b}, see p.1155. Observe also, that by (\ref{sal1})-(\ref{sal2})
\[\nu_2\ee^{hC_0}\leqslant 0.1\nu_1\ee^{-1.3997C_0}\leqslant 0.1\nu_1\ee^{-1.3997}.\]
The estimate in (\ref{lau6}) follows from  the assumptions, (\ref{lau1}), (\ref{lau2}) and Lemma \ref{Hei3aa}. Indeed,

\begin{eqnarray*}
\lefteqn{\ab{\wE (Y_j, Y_2, \cdots, Y_k)}\leqslant
2^{k-j}(\Expect\ab{Y_1}^2)^{(k-j+1)/2}\leqslant 2^{k-j}(\Expect\ab{Y_1}^2)^2(\ee^{-4.7C_0}0.009)^{(k-j-3)/2}}
\hskip 1cm\\
&\leqslant& 2^{k-j}(\ee^{2.42C_0}\nu_1^2(128\sin^4(t/2)+2\cdot2.02^2y^4)(\ee^{-4.7C_0}0.009)^{(k-j-3)/2}\\
&\leqslant&2^{k-j}\ee^{-2.58C_0}0.002\nu_1128\sin^2(t/2)(0.0091)^{k-j-3}
\\
&&+2^{k-j}8.17\nu_1^2y^2(\ee^{-2.28}0.009)^{k-j-3}\\
&\leqslant&0.16\nu_1\sin^2(t/2)(0.02)^{k-j-3}+65.4\nu_1^2y^2(0.061)^{k-j-3}.
\end{eqnarray*}
Similarly,
\begin{eqnarray*}
\ab{\Expect Y_1Y_2}&\leqslant&\ee^{0.21C_0}\Expect X_1X_2(8\sin^2(t/2)+2.02y^2)\leqslant \ee^{0.21C_0-1.5C_0}0.4\nu_1
\sin^2(t/2)\\
&&+2.5\ee^{0.21C_0}\Expect X_1X_2 y^2\leqslant 0.11011\nu_1\sin^2(t/2)+2.5\ee^{0.21C_0}\Expect X_1X_2 y^2
\end{eqnarray*}
and
\begin{eqnarray*}\Expect\ab{Y_1Y_2Y_3}&\leqslant& C_0\ee^{0.1003C_0}2.1003\Expect\ab{Y_1Y_2}\\
&\leqslant&
2.1003(\ee^{(1.1003-1.29)C_0}0.4\nu_1\sin^2(t/2)+2.5\ee^{1.4C_0}\Expect X_1X_2y^2)\\
&\leqslant&0.695\nu_1\sin^2(t/2)+5.26\ee^{1.4C_0}y^2\Expect X_1X_2.
\end{eqnarray*}
The estimates (\ref{lau4}) and (\ref{lau5}) now follow from the obvious relations $\ab{\wE(Y_1,Y_2)}\leqslant\Expect\ab{Y_1Y_2}+\Expect\ab{Y_1}\Expect\ab{Y_2}$ and $\ab{\wE(Y_1,Y_2,Y_3)}\leqslant\Expect\ab{Y_1Y_2Y_3}+\Expect\ab{Y_1Y_2}\Expect\ab{Y_3}+\Expect\ab{Y_1}\Expect\ab{Y_2Y_3}+\Expect\ab{Y_1}\Expect\ab{Y_2}
\Expect\ab{Y_3}$.
The estimates in (\ref{lau6a}) are proved similarly, with replacing $\nu_1^2$ by $\gamma$ in all steps. Also, the
estimates (\ref{lau4a}) and (\ref{lau5a}) can be verified directly. For the proof of (\ref{lau7}), observe that
\begin{eqnarray*}
\ab{(\wE (Y_j, Y_2, \cdots, Y_k))'}&\leqslant&\sum_{m=j}^k\ab{\wE(Y_j,\dots,Y'_m,\dots,Y_k)}\leqslant
\sum_{m=j}^k2^{k-j}\sqrt{\Expect\ab{Y'_m}^2}\prod_{l\ne m}^k\sqrt{\Expect\ab{Y_l}^2}\\
&\leqslant&C(k-j+1)2^{k-j}\sqrt{\nu_1}(\Expect\ab{Y_1}^2)^{3/2}(\Expect\ab{Y_1}^2)^{(k-j-3)/2}\\
&\leqslant&
C(k-j+1)2^{k-j}(\ab{\sin(t/2)}^3+y^3)(\ee^{-4.7}0.009)^{(k-j-3)/2}.
\end{eqnarray*}
$\square$
\begin{lemma}\label{HIV} Let the conditions stated in (\ref{sal1})-(\ref{sal2})  be
satisfied. Then
\[\w F_n(u)=\vfi_1(u)\vfi_2(u)\dots\vfi_n(u),\]
 where $\vfi_1(u)=\Expect e^{u X_1}$ and, for $k=2,\dots, n$,
\[\vfi_k(u)=1+ \Expect Y_k+\sum_{j=1}^{k-1}\frac{\wE(Y_j,Y_{j+1},\dots,Y_k)}
{\vfi_j(u)\vfi_{j+1}(u)\dots \vfi_{k-1}(u)}.\]
\end{lemma}
The above lemma \label{Hversion} is a part of Lemma 3.1 from \cite{H82}.

 Let $\psi_k(u)=\vfi_k(u)/(1+\nu_1(\ee^h-1))$. Then, combining Lemmas \ref{lem1} and \ref{HIV}, we get the expansion for $\psi_k(u)$.
%
\begin{lemma}\label{lemfi} Let the conditions in (\ref{sal1})-(\ref{sal2})  be
satisfied. Then, for $k=1,2,\dots,n$,
\begin{eqnarray}
\frac{1}{\ab{\vfi_k(u)}}&\leqslant&\frac{10}{9},\quad \ab{\vfi'_k}\leqslant C\label{ps1}\\
\ab{\psi_k(u)}&\leqslant& (1+20\gamma y^2)\exponent{-0.5\nu_1\sin^2(t/2)},\label{ps2}\\
\psi_k(u)&=&1+\frac{x}{n}(\ee^{\ii t}-1)+\theta C \gamma (\sin^2(t/2)+y^2),\label{ps3}\\
\psi'_k(u)&=&\ii \ee^{\ii t}\frac{x}{n}+\theta C \gamma(\ab{\sin(t/2)}+y).\label{ps3a}
\end{eqnarray}
\end{lemma}
\Proof The estimte in (\ref{ps1}) is  proved in \cite{CeVe15b}, Lemma 4.3. From Lemma \ref{HIV}, equation
(\ref{au1}) and Lemma \ref{lem1}, it follows that
\begin{eqnarray}
\psi_k(u)&=&1+\frac{x}{n}(\ee^{\ii t}-1)+\theta(0.1\nu_1\sin^2(t/2)+1.01\nu_2\ee^{0.2C_0}y^2)\nonumber\\
&&+\theta\sum_{j=1}^{k-1}\bigg(\frac{10}{9}\bigg)^{k-j}\ab{\wE(Y_j,\dots,Y_k)}\nonumber\\
&=&1+\frac{x}{n}(\ee^{\ii t}-1)+\theta (1.262\nu_1\sin^2(t/2)+20\gamma y^2).\label{apsi1}
\end{eqnarray}
Standard calculations show that
\begin{equation}
\Ab{1+\frac{x}{n}(\ee^{\ii t}-1)}^2=1-4\frac{x}{n}\bigg(1-\frac{x}{n}\bigg)\sin^2(t/2).\label{xnu}\end{equation}
From (\ref{sal1})-(\ref{sal2}), it follows that
\[\frac{9\nu_1}{10}\leqslant\frac{x}{n}\leqslant\frac{11\nu_1}{10},\qquad \frac{x}{n}\bigg(1-\frac{x}{n}\bigg)\geqslant 0.899\nu_1\sin^2(t/2), \]
so that
\[ \Ab{1+\frac{x}{n}(\ee^{\ii t}-1)}\leqslant 1-2\cdot\frac{9}{10}\nu_1\bigg(1-\frac{11}{10}\cdot 0.002\bigg).\]
Substituting this estimate into (\ref{apsi1}), we obtain
\[\ab{\psi_k(u)}\leqslant 1-2\cdot0.89\nu_1 \sin^2(t/2)+1.2\nu_1\sin^2(t/2)+20\gamma y^2\leqslant (1-0.5\nu_1\sin^2(t/2))
(1+20\gamma y^2)\]
and (\ref{ps2}) follows. The expansion in (\ref{ps3}) can be proved similarly as in Lemma \ref{lem1}, using (\ref{lau4a}), (\ref{lau5a}) and (\ref{lau6a}). For the proof of (\ref{ps3a}), observe that
\begin{eqnarray*}
\vfi'_k(u)&=&\Expect Y'_k+\frac{(\wE(Y_{k-1},Y_k))'}{\vfi_{k-1}(u)}+\frac{(\wE(Y_{k-2},Y_{k-1},Y_k))'}{\vfi_{k-2}(u)\vfi_{k-1}(u)}\\
&&+\sum_{j=1}^{k-3}\frac{(\wE(Y_j,\dots,Y_k))'}{\vfi_{j}(u)\cdots\vfi_{k-1}(u)}-\sum_{j=1}^{k-1}
\frac{\wE(Y_j,\dots,Y_k)}{\vfi_j(u)\cdots\vfi_{k-1}(u)}\sum_{m=j}^{k-1}\frac{\vfi'_m(u)}{\vfi_m(u)}.
\end{eqnarray*}
Applying to the last identity estimates from Lemma \ref{lem1} and (\ref{ps1}), we obtain

\begin{equation}
\vfi'_k(u)=\ii \ee^u\nu_1+\theta C\gamma(\ab{\sin(t/2)}+y).\label{ivfi}
\end{equation}
Dividing the last expression by $1+\nu_1(\ee^h-1)$ and applying (\ref{au1}), the identity in (\ref{ps3a}) follows.
$\square$

To prove Theorem \ref{THM2}, we need  some additional notations. Let $\ZZ$ denote the set of integers. For an integer-valued measure $M$,  define the Kolmogorov and total variation metrics  by
\[\ab{M}_K:=\sup_{m\in\ZZ}\Ab{\sum_{k=-\infty}^mM\{k\}},\qquad \norm{M}:=\sum_{k=-\infty}^\infty\ab{M\{k\}},\]
respectively. Note that $\ab{M}_K\leqslant\norm{M}$. The following inversion formula allows to estimate total variation of $M$ via its Fourier-Stieltjes transform $\w M(\ii t)=\sum_{k=-\infty}^\infty\ee^{\ii tk}M\{k\}$:


\begin{lemma}\label{varijotas} Let the  measure $M$ be concentrated on $\ZZ$ and
$\sum_{k\in\ZZ}\ab{k}\ab{M{k}}<\infty$. Then, for any
$a\in\RR$ and $b>0$, the following inequality holds:
\begin{equation*}
    \norm{M}\leqslant
(1+b\pi)^{1/2}\Biggl( \frac{1}{2\pi}\int_{-\pi}^\pi
  \biggl( \Ab{\w M(\ii t)}^2+\frac{1}{b^2}
  \Ab{\Bigl(\ee^{-\ii ta}\w M(\ii t)\Bigr)'}^2
   \biggr)\,\dd t
  \Biggr)^{1/2}.
 \label{var}\end{equation*}
\end{lemma}
Lemma \ref{varijotas} is a well known inversion inequality for lattice distributions. Its proof can be found, for example, in \cite{Ce16}, Lemma 5.1. 


 The following simple equality allows to switch from rv to its conjugate rv.

\begin{lemma}\label{GZ} Let $G$ be an integer-valued distribution and let $\w G(\ii t)$ be its characteristic function.
 Then, for any bounded $z>0$ and an integer $m$,
\begin{equation}
\ee^{zm}G\{m\}=\frac{1}{2\pi}\int_{-\pi}^\pi\w G(\ii t+z)\ee^{-\ii tm}\,\dd t.\label{conz}
\end{equation}
\end{lemma}

\Proof We introduce conjugate distribution $G_z\{m\}=\ee^{zm}G\{m\}/\w G(z)$. Observe that the characteristic function
of $G_z$ is equal to $\w G(\ii t+z)/\w G(z)$. It remains to apply formula of inversion
\[G_z\{m\}=\frac{1}{2\pi}\int_{-\pi}^\pi \frac{\w G(\ii t+z)}{\w G(z)}\ee^{-\ii tm}\,\dd t.\]
$\square$

\begin{lemma}\label{poisonx} Let $z$ be a solution of the saddle point equation $(n\nu_1(\ee^z-1))'=x$. Then
\[
\frac{(n\nu_1)^x\ee^{-n\nu_1}}{x!}=\exponent{n\nu_1(\ee^z-1)-zx}\,\frac{x^x\ee^{-x}}{x!}.\]
\end{lemma}
\Proof Observe that $n\nu_1\ee^z=x$ and
\[\exponent{n\nu_1(\ee^{\ii t+z}-1)}=\exponent{n\nu_1(\ee^z-1)}
\exponent{n\nu_1(\ee^{\ii t +z}-\ee^z)}.
\]
Therefore, by (\ref{conz}).
\[
\frac{(n\nu_1)^x\ee^{-n\nu_1}}{x!}=\exponent{n\nu_1(\ee^z-1)-zx}\frac{1}{2\pi}\int_{-\pi}^\pi
\ee^{x(\ee^{\ii t}-1)}\ee^{-\ii tx}\,\dd t.\quad \square\]


The proofs of the next two lemmas are  similar to that of Lemma \ref{poisonx}.
\begin{lemma}\label{nbx} Let $w$ be a solution of the saddle point equation $(\ln\w{\NB}(r,\qubar)(w))'=x$, $u=w+\ii t$.  Then
\[
\NB(r,\qubar)\{x\}=\ee^{-wx}\w{\NB}(r,\qubar)(w)\NB(r,r/(r+x))\{x\}
\]
and
\[\frac{\w{\NB}(r,\qubar)(u)}{\w{\NB}(r,\qubar)(w)}=\bigg(\frac{1-x/(r+x)}{1-x\ee^{\ii t}/(x+r)}\bigg)^r=\w{\NB}(r,r/(x+r))(\ii t).\]
\end{lemma}
\begin{lemma}\label{bix} Let $\tilde{h}$ be a solution of the saddle point equation $(\ln\w{\BI}(N,\tilde{p})(\tilde{h}))'=x$,
$u=\tilde{h}+\ii t$.  Then
\[
\BI(N,\tilde{p})\{x\}=\ee^{-\tilde{h}x}\w{\BI}(N,\tilde p)(\tilde{h})\BI(N,x/N)\{x\}
\]
and
\[\frac{\w{\BI}(N,\tilde{p})(u)}{\w{\BI}(N,\tilde{p})(\tilde{h})}=\bigg(1+\frac{x}{N}(\ee^{\ii t}-1)\bigg)^N=\w{\BI}(N,x/N)(\ii t).\]
\end{lemma}


We also need some estimates for gamma functions.
\begin{lemma}\label{gamaf} For all positive real numbers $x\geqslant 1$ we have
\[x^x\ee^{-x}\sqrt{2\pi(x+0.16)}<\Gamma(x+1)<x^x\ee^{-x}\sqrt{2\pi(x+0.18)}.\]
\end{lemma}
Lemma \ref{gamaf} is a slightly rougher version of Theorem 1.6 in \cite{Bat08}.

\section{Proofs}

We use the same notations as in the previous sections.

\emph{\textbf{Proof of Theorem \ref{THM1}}}. Note that $h$ is not a solution of the saddle point equation for $S_n$, which is very complicated. Instead, $h$ is a solution for the binomial distribution with the same mean. The heuristics for such a 'truncated' version of  equation are the following: the assumptions in (\ref{sal1})-(\ref{sal2}) imply  that the main probabilistic mass of $S_n$ is concentrated at 0 and 1.
From Lemma \ref{poisonx} and the equation (\ref{conz}), it follows that
\begin{eqnarray}
\Prob(S_n=x)&=&\ee^{-hx}\frac{1}{2\pi}\int_{-\pi}^\pi\prod_{j=1}^n\vfi_j(u)\ee^{-\ii tx}\,\dd t\nonumber\\
&=&
\frac{(n\nu_1)^x\ee^{-n\nu_1}}{x!}\exponent{-xh+n\ln(1+\nu_1(\ee^h-1))+zx-n\nu_1(\ee^z-1)}\nonumber\\
&&\times\bigg(1+\bigg(\frac{x^x\ee^{-x}}{x!}\bigg)^{-1}\frac{1}{2\pi}\int_{-\pi}^\pi
\Big(\prod_{k=1}^n\psi_k(u)-\ee^{x(\ee^{\ii t}-1)}\Big)\ee^{-\ii tx}\,\dd t\bigg).\label{trr}
\end{eqnarray}
Observe that  $z=\ln(1+y)$ and
\[\ee^z=\frac{x}{n\nu_1}=1+y,
\quad h=\ln(1+y)-\ln\Big(1-\frac{\nu_1y}{1-\nu_1}\Big),\quad 1+\nu_1(\ee^h-1)=\Big(1-\frac{\nu_1y}{1-\nu_1}\Big)^{-1}.\]
Therefore, $-xh+n\ln(1+\nu_1(\ee^h-1))+zx-n\nu_1(\ee^z-1)=\Lambda(y)$. Next, from the first inequality in (\ref{xnu}), the estimates in (\ref{ps2})
and (\ref{ps3}) and the trivial inequality $a\ee^{-Ca}\leqslant C$, it follows that
\begin{eqnarray}
\lefteqn{\Ab{\prod_{k=1}^n\psi_k(u)-\exponent{x(\ee^{\ii t}-1)}}}\hskip 1cm\nonumber\\
&\leqslant& (1+20\gamma y^2)^{n-1}
\ee^{-C(n-1)\nu_1\sin^2(t/2)}\sum_{k=1}^n\ab{\psi_k(u)-\ee^{(x/n)(\ee^{\ii t}-1)}}\nonumber\\
&\leqslant&C(1+20\gamma y^2)^n\ee^{-Cn\nu_1\sin^2(t/2)}Cn\bigg(\bigg(\frac{x}{n}\bigg)^2\sin^2(t/2)+\gamma(\sin^2(t/2)+y^2)\bigg)\nonumber\\
&\leqslant&C(1+20\gamma y^2)^n\ee^{-Cn\nu_1\sin^2(t/2)}Cn\gamma(\sin^2(t/2)+y^2)\nonumber\\
&\leqslant&
(1+20\gamma y^2)^n\ee^{-Cn\nu_1\sin^2(t/2)}Cn\gamma((n\nu_1)^{-1}+y^2)\nonumber\\
&=&(1+20\gamma y^2)^n\ee^{-Cn\nu_1\sin^2(t/2)}C\gamma((\nu_1)^{-1}+ny^2).
\label{delta}
 \end{eqnarray}
 One can use Stirling's formula or  Lemma \ref{gamaf} to get the estimate
 \[\bigg(\frac{x^x\ee^{-x}}{x!}\bigg)^{-1}\leqslant C\sqrt{x}\leqslant C\sqrt{n\nu_1}.\]
 Substituting all the above estimates into (\ref{trr}) and observing that
 \begin{equation}\label{sinint}\int_{-\pi}^\pi\exponent{-Cn\sin^2(t/2)}\,\dd t\leqslant C(n\nu_1)^{-1/2}\end{equation}
  the proof is completed. $\square$

\emph{\textbf{Proof of Theorem \ref{THM2}}.}
 The solution to the saddle point equation for the Poisson law with parameter $n\lambda^{*}$ should satisfy equation $(n\lambda^{*}(\ee^h-1))'=x$, that is, $n\lambda^{*}\ee^h=x$,
 where $x$ is defined in \eqref{au1}. The
  standard approach to LD requires to solve this equation with respect to $h$. Instead, we use $h$ defined by (\ref{au1}) and solve equation with respect to $\lambda^{*}$, obtaining
\begin{equation}\label{lamb}
\lambda^{*}=\nu_1\bigg(1-\frac{\nu_1y}{1-\nu_1}\bigg).
\end{equation}
As before, $y=(x-n\nu_1)/(n\nu_1)$, defined in \eqref{neqn1}. Let us denote the distribution of $Z^{*}_n$ by $\Pi^{*}$.  Recall that $F_n$ is the distribution of $S_n$. Let
\[\Pi_h^{*}\{k\}:=\frac{\ee^{hk}\Pi^{*}\{k\}}{\ee^{n\lambda^{*}(\ee^h-1)}},\qquad F_h\{k\}:=
\frac{\ee^{hk}F_n\{k\}}{(1+\nu_1(\ee^h-1))^n},\qquad k=0,1,\dots.\]
We have
\begin{eqnarray*}
\Prob(S_n\geqslant x)&=&\sum_{k\geqslant x}F_n\{k\}=\exponent{n\ln(1+\nu_1(\ee^h-1))-hx}\sum_{k\geqslant x}\ee^{-h(k-x)}F_h\{k\},\\
\Prob(Z^{*}_n\geqslant x)&=&\sum_{k\geqslant x}\Pi^{*}_n\{k\}=\exponent{n\lambda^{*}(\ee^h-1)-hx}\sum_{k\geqslant x}\ee^{-h(k-x)}\Pi^{*}_h\{k\}.
\end{eqnarray*}
Therefore,
\begin{eqnarray}
\lefteqn{\frac{\Prob(S_n\geqslant x)}{\Prob(Z_n^{*}\geqslant x)}}\hskip 0.5cm\nonumber\\
&=&\ee^{\Lambda^{*}(y)}\bigg(1+\Big(\sum_{k\geqslant x}\ee^{-h(k-x)}\Pi^{*}_h\{k\}\Big)^{-1}\sum_{k\geqslant x}\ee^{-h(k-x)}(F_h\{k\}-\Pi_h^{*}\{k\})\bigg).
\label{AA}
\end{eqnarray}
It is easy to check that $\Pi_h^{*}$ has Poisson distribution with parameter $x$, since $n\lambda^{*}\ee^h=x$. Therefore,
\begin{equation}
\sum_{k\geqslant x}\ee^{-h(k-x)}\Pi^{*}_h\{k\}\geqslant\Pi_h^{*}\{x\}=\frac{x^x\ee^{-x}}{x!}\geqslant \frac{C}{\sqrt{x}}\geqslant\frac{C}{\sqrt{n\nu_1}}.\label{AB}
\end{equation}
Abel's partial summation formula (see, for example, (1.44) in \cite{Ce16}) allows to prove
  \begin{equation}
  \Ab{\sum_{k\geqslant x}\ee^{-h(k-x)}(F_h\{k\}-\Pi_h^{*}\{k\})}\leqslant\ab{F_h-\Pi^{*}_h}_K\leqslant\norm{F_h-\Pi_h^{*}}.
  \label{AC}\end{equation}
    To estimate $\norm{F_h-\Pi_h^{*}}$, we apply Lemma \ref{varijotas}. Let $\w F_h(\ii t)$ and $\w\Pi_h^{*}(\ii t)$ denote characteristic functions of $F_h$ and $\Pi_h^{*}$, respectively. Then
  \[\w \Delta(t):=\w F_h(\ii t)-\w\Pi_h^{*}(\ii t)=\prod_{k=1}^n\psi_k(u)-\ee^{x(\ee^{\ii t}-1)}.
    \]
    The definition of $\psi_k(u)$ is given before Lemma \ref{lemfi}.
 We will estimate $\w\Delta'(t)$. To shorten our expressions, let
 \[a_k:=\psi_k(u)\exponent{-\ii tx/n},\qquad b_k:=\exponent{(x/n)(\ee^{\ii t}-1)}\exponent{-\ii tx/n}.\]
 Then
 \[\w\Delta \ee^{-\ii tx}=\prod_{k=1}^na_k-\prod_{k=1}^n b_k.\]
Using the inequalities
\begin{eqnarray*}
\ab{\sin(t/2)}\ee^{-Cn\nu_1\sin^2(t/2)}&\leqslant& C(n\nu_1)^{-1/2}\ee^{-Cn\nu_1\sin^2(t/2)},\\
\ee^{-C(n-1)\nu_1\sin^2(t/2)}&\leqslant& C \ee^{-Cn\nu_1\sin^2(t/2)},\end{eqnarray*}
and the equations (\ref{ps3a}) and (\ref{delta}),  we get, after some routine calculations,

\begin{eqnarray}
\Ab{(\w\Delta(t)\ee^{-\ii tx})'}&\leqslant&\sum_{k=1}^n\Ab{(a'_k-b'_k)\prod_{l\ne k}^na_l}+\sum_{k=1}^n\Ab{b'_k\Big(
\prod_{l\ne k}^n a_k-\prod_{l\ne k}^n b_k\Big)}\nonumber
\end{eqnarray}
\begin{eqnarray}
&\leqslant&C(1+20\gamma y^2)\ee^{-Cn\nu_1\sin^2(t/2)}\Big(\sum_{k=1}^n\ab{a'_k-b'_k}+n\cdot\frac{x}{n}\ab{\sin(t/2)}\gamma
(\nu_1^{-1}+ny^2)\Big)\nonumber\\
&\leqslant&C(1+20\gamma y^2)\ee^{-Cn\nu_1\sin^2(t/2)}\bigg(n\gamma\bigg(\frac{1}{\sqrt{n\nu_1}}+y\bigg)+\sqrt{n\nu_1}(\nu_1^{-1}+ny^2)
\bigg).\label{zz}
\end{eqnarray}

 Observe, that $\ab{\w \Delta(t)}$ is already estimated in (\ref{delta}). Therefore applying Lemma \ref{varijotas} with $a=x$ and $b=\max(1,\sqrt{n\nu_1})$ and noting that $\sqrt{n}y/\sqrt{\nu_1}\leqslant \nu_1^{-1}+ny^2$, we obtain
 \[\norm{F_h-\Pi_h^{*}}\leqslant C(1+20\gamma y^2)^n\gamma(\nu_1^{-1}+ny^2).\]
 Substituting the last estimate and the estimate (\ref{AB}) into (\ref{AA}), we complete the proof. $\square$

\emph{\textbf{Proof of Theorem \ref{THM3}.}}  The main step in the proof, allowing simplification of estimates, is the following. When we deal with the choice of parameter for conjugate distribution of $S_\xi$, we do not try to solve the corresponding saddle point equation. Instead, we use solution $w$ to the saddle point equation of  the negative binomial approximation, that is,
 $(\ln\w{\NB}(r,\qubar)(w))'=x$.
 From the assumption of the theorem, it follows that
\begin{equation}
p=o(1),\quad x=np^2(1+o(1)),\quad \pbar=2p(1+o(1)),\quad r=\frac{np}{2}(1+o(1)).
\label{neig1}
\end{equation}
Saddle point solution $w$ satisfies $\pbar\ee^w=x/(r+x)$. Therefore,
\begin{equation}
\ee^w-1=\frac{x-np^2}{np^2}(1+o(1))=o(1),\quad \ab{\ee^u-1}\leqslant 2\ab{\sin(t/2)}+(\ee^w-1)\leqslant 2+o(1).
\label{neig2}
\end{equation}
For 2-runs, we have an explicit version of Heinrich's lemma, see \cite{PeCe10} Lemmas 1 and 2 (and beware of a missprint). Let $\Psi_\xi(u):=\Expect \exponent{uS_\xi}$. Then,
\begin{eqnarray}
\Psi_\xi(u)&=&\prod_{k=1}^n f_k(u),\quad f_1(u)=1+p^2(\ee^u-1),\label{neig3}\\
f_k(u)&=&1+p^2(\ee^{u}-1)+\sum_{j=1}^{k-1}\frac{(\ee^u-1)^{k-j+1}p^{k-j+2}(1-p)^{k-j}}{f_j(u)\cdots f_{k-1}(u)}.\label{neig4}
\end{eqnarray}
Next steps in the proof are very similar to the steps from the proof of Lemma 7 in \cite{PeCe10} and, therefore, will not be discussed in detail. By induction, we prove that $\ab{f_k(u)-1}\leqslant Cp^2\ab{\ee^u-1}^2$, then we prove $\ab{f_k(u)-1-p^2(\ee^{u}-1)}\leqslant Cp^3\ab{\ee^u-1}^3$, and  finally obtain
\[
\ln\Psi_\xi(u)=np^2(\ee^u-1)+\frac{np^3(2-3p)-2p^3(1-p)}{2}\,(\ee^u-1)^2+\theta Cnp^4\ab{\ee^u-1}^3.\]
Similarly,
\[\w{\NB}(r,\qubar)(u)=\Exponent{r\bigg(\frac{\qubar}{\pbar}\bigg)(\ee^u-1)+r\bigg(\frac{\qubar}{\pbar}\bigg)^2(\ee^u-1)^2+\theta Cr
\bigg(\frac{\qubar}{\pbar}\bigg)^3\ab{\ee^u-1}^3}.\]
Taking into account (\ref{neig1}), (\ref{neig2}) and the  definitions of $r,\pbar,\qubar$ defined in \eqref{NBpar}, we get
\begin{equation}
\Psi_\xi(u)/\w{\NB}(r,\qubar)(u)=\exponent{Cnp^4\ab{\ee^u-1}^3}=\exponent{Cnp^4\ab{\sin(t/2)}^3+o(1)}.
\label{neig5}
\end{equation}
Observe also that
\begin{eqnarray}
\ab{\w{\NB}(r,x/(r+x))(\ii t)}&=&\Ab{\Exponent{r\sum_{j=1}^\infty\bigg(\frac{x}{r+x}\bigg)^j\frac{(\ee^{\ii tj}-1)}{j}}}\nonumber\\
&\leqslant&
\Exponent{-\frac{2rx}{r+x}\sin^2(t/2)}\leqslant \exponent{-np^2\sin^2(t/2)},
\label{neig6}
\end{eqnarray}
if $n$ is sufficiently large. Applying Lemmas \ref{GZ} and \ref{nbx}, we obtain
\begin{eqnarray}
\lefteqn{\Prob(S_\xi=x)=\frac{\ee^{-wx}}{2\pi}\int_{-\pi}^\pi\Psi_\xi(u)\ee^{-\ii tx}\,\dd x}\hskip 0.5cm\nonumber\\
&=&\NB(r,\qubar)\{x\}\bigg(
\ee^{-wx}\w{\NB}(r,\qubar)(w)\NB\Big(r,\frac{r}{r+x}\Big)\{x\}\bigg)^{-1}\frac{\ee^{-wx}}{2\pi}\int_{-\pi}^\pi\Psi_\xi(u)\ee^{-\ii tx}\,\dd x
\nonumber\\
&=&\NB(r,\qubar)\{x\}\bigg(\NB\Big(r,\frac{r}{r+x}\Big)\{x\}\bigg)^{-1}
\frac{1}{2\pi}\int_{-\pi}^\pi\frac{\w{\NB}(r,\qubar)(u)}{\w{\NB}(r,\qubar)(w)}
\frac{\Psi_\xi(u)}{\w{\NB}(r,\qubar)(u)}\,\ee^{-\ii t x}\,\dd x\nonumber\\
&=&\NB(r,\qubar)\{x\}\bigg(1+\bigg(\NB\Big(r,\frac{r}{r+x}\Big)\{x\}\bigg)^{-1}\nonumber\\
&&\times\frac{1}{2\pi}
\int_{-\pi}^\pi\w{\NB}\Big(r,\frac{r}{r+x}\Big)(\ii t)\bigg(\frac{\Psi_\xi(u)}{\w{\NB}(r,\qubar)(u)}
-1\bigg)\ee^{-\ii t x}\,\dd x\bigg)\label{neigneig}
\end{eqnarray}
 It follows from (\ref{neig5}) and (\ref{neig6}) that, for sufficiently large $n$,
\begin{eqnarray*}
\lefteqn{\int_{-\pi}^\pi\Ab{\w{\NB}\Big(r,\frac{r}{r+x}\Big)(\ii t)\bigg(\frac{\Psi_\xi(u)}{\w{\NB}(r,\qubar)(u)}
-1\bigg)}\dd t}\hskip 0.3cm\\
&\leqslant&
C\int_{-\pi}^\pi\ee^{-2np^2\sin^2(t/2)+Cnp^4\ab{\sin(t/2)}^3}(np^4\ab{\sin(t/2)}^3+o(1))\dd t\\
&\leqslant&C\int_{-\pi}^\pi\ee^{-1.5np^2\sin^2(t/2)}(np^4\ab{\sin(t/2)}^3+o(1))\dd t\\
&\leqslant&C\int_{-\pi}^\pi
\ee^{-np^2\sin^2(t/2)}(np^4(np^2)^{-3/2}+o(1))\dd t=\frac{o(1)}{\sqrt{np^2}}.
\end{eqnarray*}
Using Lemma \ref{gamaf} and the equation (\ref{neig1}), we get
\[
\bigg(\NB\Big(r,\frac{r}{r+x}\Big)\{x\}\bigg)^{-1}\leqslant C\sqrt{\frac{x}{1+x/r}}\leqslant C\sqrt{np^2}.
\]
Substituting the last two estimates into (\ref{neigneig}),  the proof of Theorem \ref{THM3} follows. $\square$

\emph{\textbf{Proof of Theorem \ref{THM4}}.} It is easy to check that for $\tilde S$ the following version of Lemma \ref{HIV}
can be used:
\begin{eqnarray}\tilde\Psi(u)&:=&\Expect\exponent{u\tilde S}=\prod_{k=1}^n g_k(u),\quad g_1(u)=1+\alpha(\ee^u-1),\nonumber\\
g_k(u)&=&1+\alpha(\ee^u-1)+\sum_{j=1}^{k-1}\frac{(-1)^{k-j}\alpha^{k-j+1}(\ee^u-1)^{k-j+1}}{g_j(u)\cdots g_k(u)},\quad k\geqslant 2. \label{bin1}
\end{eqnarray}
 Here, $u=\tilde{h}+\ii t$ and $\tilde{h}$ is now a solution of the saddle point equation for the binomial distribution:
\[\frac{N\tilde p\ee^{\tilde{h}}}{1+\tilde{p}(\ee^{\tilde h}-1)}=x.\]
For sufficiently large $n$, we get
\[\ee^{\tilde h}-1=\frac{x-n\alpha}{n\alpha}(1+o(1)),\quad x=n\alpha(1+o(1)),\quad N=\frac{n}{3}(1+o(1)),\quad\tilde{p}=3\alpha(1+o(1)).\]
Recursive formula (\ref{bin1}) allows us to prove
\[\tilde\Psi(u)=\Exponent{n\alpha(\ee^u-1)-\frac{\alpha^2(3n-2)}{2}(\ee^u-1)^2+C\theta n\alpha^3\ab{\sin(t/2)}^3+o(1)}.\]
 Similar expression holds for the binomial distribution
 \[\w{\BI}(N,\tilde{p})(u)=\Exponent{N\tilde{p}(\ee^u-1)-\frac{N\tilde{p}^2}{2}(\ee^u-1)^2+\theta C N\tilde{p}^3\ab{\ee^u-1}^3}.\]
The choice of parameters ensures that $N\tilde{p} =n\alpha$ and
\[
\ab{N\tilde{p}^2-(3n-2)\alpha^2}=\alpha^2\frac{(3n-2)\epsilon/n}{n/(3n-2)+\epsilon/n}\leqslant C\alpha^2.\]
Applying Lemma \ref{gamaf}, we have for sufficiently large $n$,
\[\BI(N,N/x)\{x\}\geqslant C\sqrt{x}\geqslant C\sqrt{n\alpha}.\]
The proof of Theorem \ref{THM4} can be completed by writing the analogue of (\ref{neigneig}) and using similar arguments as in the proof Theorem \ref{THM3}. $\square$

\vskip 0.5cm
\section*{Acknowledgment}
The main part of this paper was written during the first
author's stay at the Department of Mathematics, IIT Bombay, during
January-February, 2018. The first author would like to thank the members
of the Department for their hospitality.
We are grateful to the referees for useful remarks, which helped to improve the paper.

\end{document}